\def\be{\begin{equation}}
\def\ee{\end{equation}}
\def\ba{\begin{array}}
\def\ea{\end{array}}

\def\l{\lambda}
\def\p{\partial}

\def\Rb{{I\!\! R}}

\def\Cb{\ \hbox{\vrule width 0.6pt height 7pt depth 0pt
		      \hskip -3.5 pt} C}
\documentstyle[12pt]{article}
\begin{document}
\parskip=4pt
\parindent=18pt
\baselineskip=22pt
\setcounter{page}{1}
\centerline{\Large\bf Symplectic Geometry on Quantum Plane}
\vspace{6ex}
\begin{center}
{\large  Sergio Albeverio}\footnote{SFB  237; BiBoS; CERFIM (Locarno); 
Acc.Arch., USI (Mendrisio)} and
{\large  Shao-Ming Fei}\footnote{Institute of Physics,
Chinese Academy of Science, Beijing.}
\end{center}
\begin{center}
Institut f\"ur Angewandte Mathematik, Universit\"at Bonn, D-53115 Bonn\\
and\\
 Fakult\"at f\"ur Mathematik, Ruhr-Universit\"at Bochum, D-4478 Bochum
\end{center}
\vskip 1 true cm
\parindent=18pt
\parskip=6pt
\begin{center}
\begin{minipage}{5in}
\vspace{3ex}
\centerline{\large Abstract}
\vspace{4ex}
A study of symplectic forms associated with two dimensional quantum
planes and the quantum sphere in a three dimensional orthogonal quantum
plane is provided. The associated Hamiltonian
vector fields and Poissonian algebraic relations are made explicit.
\end{minipage}
\end{center}
\newpage

Originated from the investigation of trigonometric and hyperbolic 
solutions of the quantum Yang-Baxter equation (QYBE)\cite{yang,baxter}, the
quantum groups \cite{qg} have attracted much attention to
mathematicians and physicists. This is because of their relations with
many physical aspects such as exactly soluble models
in statistical mechanics, conformal field theory, integrable model field
theories \cite{baxter,app} and the fact that they have a rich mathematical
structure, see e.g., \cite{qbook}. In addition
by defining a consistent differential calculus on the non-commutative
space of the quantum groups \cite{woro}, quantum 
groups supply concrete examples of non-commutative differential 
geometry \cite{ngeom}.

The quantum plane is a space upon
which the quantum group acts by linear transformations and
whose coordinates belong to a non-commutative associative algebra 
\cite{manin}. By using the differential geometry method in 
\cite{woro} and interpreting the dual space of the quantum plane 
as differentials of the coordinates, covariant differential
calculus on the quantum plane has been developed, so that
the quantum plane provides
a simple example for non-commutative differential geometry \cite{wess}.
From the differential calculus on the three dimensional quantum plane, the
quantum Schr\"odinger equation and the $q$-deformation of the hydrogen
atom have been studied \cite{song}.
$q$-deformed integrals on the quantum plane have also been studied,
together with quantum versions of Cauchy's and Stokes' theorems \cite{qint}.

In this paper we study symplectic geometry, Poisson algebraic
structures and dynamics on quantum planes.
We give explicit expressions of the symplectic forms associated with 
Hamiltonian vector fields and the corresponding Poisson 
algebraic relations for the two dimensional
quantum plane and the quantum sphere in a three dimensional orthogonal
quantum plane. It is found that 
for some quantum planes the symplectic structures exist for
a special value of the deformation parameter $q$.

The quantum plane is defined in terms of $n$ algebraic variables $x^i$,
$i=1,2,...,n$, building a noncommutative algebra $A$ over $\Cb$,
satisfying the commutation relations \cite{manin,wess},
\be\label{1}
x^ix^j=qx^jx^i,~~~~~i<j
\ee
where $q$ is a complex number, $q\neq 0$. In addition one has algebraic
variables $\xi^i$, $i=1,2,...,n$, building a noncommutative algebra
$A^\p$, satisfying
\be\label{2}
\xi^i\xi^i=0,~~~\xi^i\xi^j=-\frac{1}{q}\xi^j\xi^i,~~~~i<j.
\ee
By definition the space $GL_q(n)$ of ``quantum matrices" consists of
$x\times n$ matrices $M=(m_{ij})_{i,j=1,...,n}$ with $m_{ij}$
(non-commutative) algebraic elements, acting on $(x^1,...x^n)$ resp.
$(\xi^1,...\xi^n)$ commutating with the $x^k$ resp. $\xi^k$ and
preserving the commutation relations (\ref{1}) resp. (\ref{2}).
In \cite{wess} an interpretation of (\ref{1}), (\ref{2}) is given in
terms of a ``differential calculus", covariant with respect to $GL_q(n)$.

More generally one can consider variables $x^i$, $i=1,...,n$, belonging
to a non-commutative algebra ${\cal A}$, satisfying the commutation
relations (with the usual summation convention):
\be\label{4}
x^ix^j-B^{ij}_{kl}x^kx^l=0,
\ee
or (in tensor product notation), $(E_{12}-B_{12})x_1x_2=0$, 
where $E$ is the $n^2\times n^2$ unit matrix,
$B$ is a $n^2\times n^2$-matrix with entries in $\Cb$,
$x_1$ and $x_2$ are two representations of the ``column vector"
$\left(\begin{array}{l}x^1\\ \vdots\\ x^n\end{array}\right)$.

Let ${\cal F}$ denote the ring of all
functions of the $x^i$, $i=1,...,n$, $x^i\in {\cal A}$,
defined by formal power series in the $x^i$, with
coefficients in $\Cb$. Let ${\cal A}^\partial$ be a non commutative
algebra, generated by $\xi^i\equiv dx^i$, $i=1,..,n$, and we suppose
that ${\cal A}^\partial$ operates on ${\cal F}$ in the sense that the
operation $f\to\xi^if$ is well defined for any $f\in {\cal F}$.
Let $d$ be the ``exterior differential" on ${\cal F}$
defined by $d=\xi^i\partial_i$, $\partial_i$ being the derivative
with respect to $x^i$ (defined on the ring ${\cal F}$).
For a consistent definition of ``differential calculus" on the ``B-plane"
(\ref{4}), one needs all the commutation relations among $x^i$, $\xi^i$
and $\partial_i$. In general $x^i$, $\xi^i$ and $\partial_i$ are not
commuting with each other but are required to
have the following commutation relations:
\be\label{6}
x_1\xi_2=C_{12}\xi_1x_2,~~~~(E_{12}+D_{12})\xi_1\xi_2=0,~~~~
\partial_2\partial_1-\partial_2\partial_1F_{12}=0,
\ee
where $C$, $D$ and $F$ are matrices with entries in $\Cb$
to be determined. By requiring that $d$ be nilpotent $d^2=0$
and satisfy the Leibniz rule
$d(fg)=(df)g+f(dg)$, $f,g\in{\cal F}$, one has, as in \cite{wess}:
\be\label{11}
\partial_k x^i=\delta_k^i+C^{ij}_{kl}x^l\partial_j,~~~~
\partial_k \xi^i=D^{ij}_{kl}\xi^l\partial_j
\ee
and the following conditions for the matrices $D$, $C$, $B$ and $F$:
\be\label{10}
\ba{l}
(E_{12}-B_{12})(E_{12}-C_{12})=0,\\[3mm]
(E_{12}-B_{12})C_{23}C_{12}=C_{23}C_{12}(E_{23}-B_{23}),\\[3mm]
D_{23}C_{12}C_{23}=C_{12}C_{23}D_{12},~~~~D=C^{-1},\\[3mm]
(E_{12}-F_{12})(E_{12}-C_{12})=0,\\[3mm]
(E_{23}-F_{23})C_{12}C_{23}=C_{12}C_{23}(E_{12}-F_{12}).
\ea
\ee

The solutions of (\ref{10}) are given by the 
$\check{R}$-matrices satisfying the Yang-Baxter equation,
\be\label{ybe}
\check{R}_{12}\check{R}_{23}\check{R}_{12}
=\check{R}_{23}\check{R}_{12}\check{R}_{23},
\ee
where $\check{R}_{ij}$ denotes the matrix on the complex vector space 
$v\otimes v\otimes v$, acting as $\check{R}$ on the $i$-th and the 
$j$-th components and as the identity on the other components, i.e.
$\check{R}_{12}=\check{R}\otimes E$, $\check{R}_{23}=E\otimes\check{R}$.
The $A_n$-type $\check{R}$ matrices satisfy
\be\label{an}
(\check{R}-\l_1)(\check{R}-\l_2)=0,
\ee
with two different eigenvalues $\l_1=-q^{-1}$ and $\l_2=q$ \cite{qbook}.
The solutions of $B$, $C$, $D$ and $F$ in the system
(\ref{10}) associated with an $A_n$-type matrix are given by \cite{wess}:
\be\label{14}
B=q^{-1}\check{R},~~~C=q\check{R},~~~F=B,~~~D=C^{-1}.
\ee

The $\check{R}$ matrices satisfying the braid Yang-Baxter equation have
three different eigenvalues $\lambda_i$, $i=1,2,3$,
associated with $B_n$,
$C_n$ and $D_n$ type solutions \cite{qbook,song12}, and one has:
\be\label{bn}
(\check{R}-\l_1)(\check{R}-\l_2)(\check{R}-\l_3)=0.
\ee
For the $B_n$ type, one gets
$\l_0=q^{-2n}$, $\l_1=-q^{-1}$, $\l_2=q$.
The solutions of $B$, $C$, $D$ and $F$ in the system of equations
(\ref{10}) related to $B_n$-type $R$-matrices are given by \cite{song}
\be\label{a3}
C=q\check{R},~~~D=C^{-1},~~~ B=F=E-Q,
\ee
where $Q=(\check{R}-\l_0)(\check{R}-\l_2)(\l_1-\l_0)^{-1}(\l_1-\l_2)^{-1}$,
$\l_1\neq\l_0$, $\l_1\neq\l_2$, is one of the three projection operators
associated with $\check{R}$.

To investigate the symplectic structures of the B-quantum plane (\ref{4})
we have to define a q-deformed symplectic
form, an exterior product, an inner product and a Hamiltonian
vector field in analogy with the corresponding objects of the theory of
ordinary symplectic manifolds \cite{geom}.
Let $V_x$ (resp. $V^{\ast}_x$) denote the vector space spanned by the basis
$\{\partial_i\equiv\frac{\partial}{\partial x^i}\}$
(resp. $\{dx^i\}$), $i=1,...,n$, at $x\in{\cal A}$, so
that, with $<,>$ being the inner product between $V_x$ and $V_x^\ast$:
\be\label{16}
<\partial_i,dx^j>=\delta_i^j,
\ee
and for $f\in {\cal F}$,
\be\label{17}
<f\partial_i,x^kdx^j>=fC^{kj}_{im}x^m,
\ee
where $C$ is given by (\ref{14}) (resp. (\ref{a3})) for the $A_n$ (resp.
$B_n$)-type case. Generally a vector in $V_x$ (resp. $V_x^\ast$) has the
form $\sum_{i=1}^n a_i(x)\partial_i$, $a_i(x)\in {\cal F}$
(resp. $\sum_{i=1}^n b_i(x)dx^i$, $b_i(x)\in {\cal F}$).
The general inner product of the form $<f\partial_i,g dx^j>$, $f,g\in
{\cal F}$, can be deduced from (\ref{17}), by using linearity.

Set $V=\cup_x V_x$ (resp. $V^\ast=\cup_x V_x^\ast$).
We call $f$ (resp. $X^{\ast}$) a zero form (resp. one form) if $f\in
{\cal F}$ (resp. $X^{\ast}\in V^{\ast}$).
Let $dx^i\otimes dx^j$ be an element in
the tensor space of $V^{\ast}_x\otimes V^{\ast}_x$. Let
$\Gamma=(\Gamma_{lm}^{jk})$, $j,k,l,m=1,...,n$, with entries in $\Cb$,
be a solution of the Yang-Baxter equation (\ref{ybe}).
We define
\be\label{20}
dx^i\wedge dx^j=dx^i\otimes dx^j-\Gamma^{ij}_{kl}dx^k\otimes dx^l.
\ee
Since $\Gamma$ is a solution of QYBE,
the exterior algebra so defined is associative \cite{qwedge}.
From the definition it follows
\be\label{18}
<\partial_i, dx^j\wedge dx^k>=\delta_i^j dx^k-\delta_i^l\Gamma_{lm}^{jk}
dx^m.
\ee
Using (\ref{17}) we then have 
\be\label{19}
<\partial_i,x^jdx^k\wedge dx^l>=C_{in}^{jk}x^ndx^l-
C^{js}_{in}\Gamma^{kl}_{st}x^ndx^t.
\ee

From commutation relations as in (\ref{6}), identifying
$\xi^i\xi^j$ with $dx^i\wedge dx^j$, and the definition of wedge product
(\ref{20}), we have that the matrix $\Gamma$ must satisfy the following
relation:
\be\label{da}
(D+E)(E-\Gamma)=0.
\ee
A matrix $\Gamma$ satisfying QYBE and (\ref{da}) defines the exterior
algebra on a quantum plane. With respect to the $A_n$ case, from equation
(\ref{an}) and (\ref{da}), we have $\Gamma=\check{R}/q$ with $\check{R}$
as in (\ref{an}). Nevertheless, as $\check{R}$ satisfies the cubic relation
(\ref{bn}) for the case associated with $B_n$ type algebras, a solution
of $\Gamma$ in terms of $\check{R}$ is not obvious and not always possible.

Let $M$ denote the quantum space defined by (\ref{4}).
We call a two form $\omega$ on $M$ closed if it satisfies
\be\label{21}
d\omega=0.
\ee
Let $\rfloor$ denote the left inner product defined by 
defined by $(X\rfloor\omega)(Y)=\omega(X,Y)$
for any two vectors $X$, $Y$ and two form $\omega$ on $M$.
For any vector $X\in V$, if
\be\label{22}
X\rfloor\omega=0~~\Rightarrow~~ X=0,
\ee
we call the two form $\omega$ non-degenerate.
Condition (\ref{22}) means that
if $\omega(X,Y)=0$ for all $Y\in V$, then $X=0$.
We call a non-degenerate closed two form a symplectic form on $M$.

Let $\omega$ be the symplectic form on $M$. For $X\in V$, if
\begin{equation}\label{23}
X\rfloor\omega=-d\,f
\end{equation}
for some $f\in {\cal F}$, we call $X\equiv X_f$ the Hamiltonian vector
field associated with $f$.

Let $X_f$, $X_g$ be the Hamiltonian vector fields associated with $f$
and $g$ respectively, $f,g\in {\cal F}$. We define the Poisson bracket
of $f$ and $g$ by
\be\label{24}
[f,g]=-X_f g.
\ee

For a given Hamiltonian $H\in {\cal F}(M)$, the corresponding dynamics
is given by the following equation of motions:
\be\label{25a}
\frac{dx^i}{dt}\equiv \dot{x}^i=[x^i, H].
\ee

{\sf Remark}: The formula (\ref{24}) is formally the same as the one in
symplectic geometry theory on usual manifolds or supermanifolds with
$U$-numbers \cite{super}. Nevertheless, in the quantum plane case,
generally relations like $[f,g]=\pm [g,f]$ for $f,g\in{\cal F}$ do not
hold.

We first investigate the symplectic structure on two dimensional quantum
plane. In this case the matrix $\check{R}$ is given by
\be\label{25}
\check{R}=\left(
\ba{cccc}
q&0&0&0\\[2mm]
0&q-q^{-1}&1&0\\[2mm]
0&1&0&0\\[2mm]
0&0&0&q
\ea\right).
\ee
By (\ref{6}), (\ref{11}) and (\ref{14}) we have all the
commutation relations such as, with coordinates $x$ and $y$,
\be\label{27}
xy=qyx,
\ee
\be\label{28}
\xi^2=\eta^2=0,~~~\xi\eta=-\frac{1}{q}\eta\xi,
\ee
where $\xi=dx$, $\eta=dy$.

The symplectic form is
\be\label{29}
\omega=\xi\wedge\eta.
\ee
Explicitly from (\ref{20}) we have
\be\label{30}
\omega=\xi\wedge\eta=q^{-2}\xi\otimes\eta-q^{-1}
\eta\otimes\xi=-\eta\wedge\xi/q,
\ee
where $\Gamma=\check{R}/q$ is used.

By formula (\ref{18}) and (\ref{23}) we have the Hamiltonian vectors
$X_x$ (resp. $X_y$) associated with $x$ and (resp. $y$),
\be\label{31}
X_x=q\partial_y,~~~X_y=-q^2\partial_x.
\ee
Therefore the Poisson bracket of $x$ and $y$ is
\be\label{32}
[x,y]=-[y,x]/q=-q.
\ee
The dynamics on the quantum plane can be investigated using
formula (\ref{25a}).

We now consider the
quantum sphere in a three dimensional orthogonal quantum plane.
The corresponding $\check{R}$-matrix is
\be\label{36}
\check{R}=\begin{array}{c}++ \\+0 \\ +- \\ 0+ \\ 00 \\ 0- \\ -+ \\ -0 \\ --
\end{array}\stackrel{\begin{array}{ccccccccc}
++~ & +0 &~~~~ +-~~~~ & ~0+ &~~~~00 ~~~&~~0- & -+ & -0 & --
\end{array}}{\left( \begin{array}{ccccccccc}
~~~~q &  &  &  &  &  &  &  &  \\  & ~~~d~ &  &~~1~~&  &  &  &  &  \\  
&  & d(1-q^{-1}) &  & -dq^{-1/2} &  & q^{-1} &  &  
\\  & ~~~1~ &  & 0 &  &  &  &  &  \\  
&  & -dq^{-1/2} &  & ~~1~~ &  &  &  &  \\  
&  &  &  &  & d &  & ~~1~~ &  \\  &  & q^{-1} &  &  &  & 0 &  &  \\  
&  &  &  &  & ~~1~~ &  & 0 &  \\  &  &  &  &  &  &  &  & ~~q~~
\end{array}
\right)},
\ee
where $d=q-q^{-1}$ and the blank spaces mean that the corresponding
entries are zeros.
Let $x^+$, $x^0$ and $x^-$ denote the coordinates, $\xi^+=dx^+$,
$\xi^0=dx^0$ and $\xi^-=dx^-$ (resp.
$\partial_+=\frac{\partial}{\partial x^+}$,
$\partial_0=\frac{\partial}{\partial x^0}$ and
$\partial_-=\frac{\partial}{\partial x^-}$) denote the differentials
(resp. derivatives). From formulae (\ref{4}) and (\ref{a3})
one has the relations:
\be\label{37}
x^+x^0=q x^0x^+,~~~x^0x^-=q x^-x^0,~~~x^+x^--x^-x^+=(q^{-1/2}-
q^{1/2})x^0x^0.
\ee
\be\label{38}
\ba{l}
x^+\xi^+=q^2 \xi^+x^+,~~~~x^+\xi^0=q \xi^0x^++(q^2-1)\xi^+x^0,\\[3mm]
x^+\xi^-=\xi^-x^+ +(q^{-1}-q)q^{1/2}\xi^0x^0-(q^{-1}-q)(q-1)\xi^+x^-,\\[3mm]
x^0\xi^+=q\xi^+x^0,~~~~x^0\xi^0=q \xi^0x^0 +(q^{-1}-q)q^{1/2}\xi^+x^-,\\[3mm]
x^0\xi^-=q\xi^-x^0 +(q^{2}-1)\xi^0x^-,\\[3mm]
x^-\xi^+=\xi^+ x^-,~~~x^-\xi^0=q\xi^0x^-,~~~x^-\xi^-=q^2\xi^-x^-.
\ea
\ee
\be\label{39}
\ba{l}
\p_+x^+=1-(q^{-1}-q)(q-1)x^-\p_-+(q^2-1)x^0\p_0+q^2x^+\p_+\\[3mm]
\p_+x^0=(q^{-1/2}-q^{3/2})x^-\p_0+qx^0\p_+,~~~\p_+x^-=x^-\p_+\\[3mm]
\p_0x^+=(q^{-1/2}-q^{3/2})x^0\p_-+ qx^+\p_0,~~~~\p_0x^0=1+(q^2-1)x^-\p_-
+qx^0\p_0,\\[3mm]
\p_0x^-=qx^-\p_0,~~~\p_-x^+=x^+\p_-,~~~\p_-x^0=qx^0\p_-,~~~\p_-x^-
=1+q^2x^-\p_-.
\ea
\ee

The quantum sphere on plane (\ref{37}) is defined by \cite{podles},
\be\label{40}
r^2=q^{-1/2}x^+x^-+x^0x^0+q^{1/2}x^-x^+,
\ee
where $r^2\in\Cb$ is the center of the algebra (\ref{37}),
$r^2x^i=x^ir^2$, $i=+,0,-$. By using relation (\ref{38}) we have
\be\label{41}
x^0dx^0+\sqrt{q}x^-dx^++\frac{1}{\sqrt{q}}x^+dx^-=0.
\ee

When $q=1$ the relations (\ref{37}), (\ref{38}) and (\ref{39}) become
the ones of among the usual coordinates, differentials and derivatives,
while the so defined quantum sphere (\ref{40}) becomes the two
dimensional sphere $S^2$ in $\Rb^3$. Therefore (\ref{40}) is a kind of
quantum deformation of $S^2$, which results in the non-commmutativity of
coordinates. Another kind of deformation is the deformation from $S^2$ to
$S_q^2$ in $\Rb^3$. The Poisson algebra on $S^2_q$ is the $q$-deformed
Lie algebra $su_q(2)$ \cite{fg}. The latter deformation gives no problem
concerning the existence of a corresponding symplectic structure for all
values of $q$, as the Poisson algebraic structures on general two
dimensional manifolds in $\Rb^3$ can be explicitly given \cite{fa}. It
might be interesting that similar
to the fact that some usual manifolds, such as
the three dimensional sphere, admit no symplectic structures, some 
quantum deformations like (\ref{40}) may also have no symplectic
structures for some values of $q$. We now give the symplectic structure
on (\ref{40}) for a special value of $q$.

From (\ref{a3}) we have
\be\label{42}
D=\begin{array}{c}++ \\+0 \\ +- \\ 0+ \\ 00 \\ 0- \\ -+ \\ -0 \\ --
\end{array}\stackrel{\begin{array}{ccccccccc}
++~ & +0 ~~& +-~ & ~0+~ &~~00 ~~~&~~0- ~~&~~ -+ ~~& ~-0~~ &~~ --
\end{array}}{\left( \begin{array}{ccccccccc}
~~~q^{-2} &  &  &  &  &  &  &  &  \\  & 0 &  &~~q^{-1}~~&  &  &  &  &  \\  
&  & ~0~ &  &  &  & 1 &  &  \\  & ~~~q^{-1}~ &  & -d/q &  &  &  &  &  \\  
&  &  &  & ~~q^{-1}~~ &  & d/\sqrt{q} &  &  \\  
&  &  &  &  & 0 &  & ~~q^{-1}~~ &  \\  &  & 1 &  & d/\sqrt{q}  &  & a &  &  \\
&  &  &  &  & ~~q^{-1}~~ &  & -d/q &  \\  &  &  &  &  &  &  &  & ~~q^{-2}~~
\end{array}
\right)},
\ee
where $a=d(1-q^{-1})$.

Besides the case $q=1$, a solution of $\Gamma$ satisfying QYBE and
(\ref{da}) exists when $q=-1$. In this case we have $D^2=I$ and
$\Gamma=D=\check{R}^{-1}$.
After some calculations we obtain the symplectic form,
\be\label{w}
\ba{rcl}
\omega&=&\frac{1}{ir^2}(- x^+ dx^-\wedge dx^0 + x^- dx^0\wedge dx^+ 
- x^0 dx^+\wedge dx^-)\\[3mm]
&=&\frac{1}{ir^2}(- x^+ (dx^-\otimes dx^0 + dx^0\otimes dx^-)
+ x^- (dx^0\otimes dx^+ + dx^+\otimes dx^0)\\[3mm]
&&- x^0 (dx^+\otimes dx^- - dx^-\otimes dx^+)).
\ea
\ee
One can check directly that $d\omega=0$.
The corresponding vector fields associated with $x^{\pm,0}$ are:
\be\label{47}
X_{x^\pm}=-\left(x^\pm\partial_0 \pm i x^0\partial_\mp\right),
~~~~~~X_{x^0}=-\left(x^-\partial_- + x^+\partial_+\right).
\ee
By using relations (\ref{37}), (\ref{39}-\ref{41}) it is
straightforward to verify that (\ref{w}) and (\ref{47}) satisfy
equation (\ref{23}).

From formula (\ref{24}) we have the Poisson relations among the $x^{\pm,0}$:
\be\label{50}
[x^\pm,x^\mp]=\pm i x^0 ,~~~~[x^0,x^\pm]=x^\pm,~~~~[x^\pm,x^0]=x^\pm.
\ee

We have studied the symplectic geometry, Poisson algebraic
structures and dynamics on quantum planes. An interesting observation
is that for some quantum planes,
the symplectic structures may only exist for some special values of $q$.
To investigate the general relations between Poisson
algebraic structures and quantum planes, the geometrical quantizations,
like in the case of usual manifolds \cite{fa}, would be a
challenging problem. It may be expected that, similar to the case 
of the quantum algebras discussed in \cite{fg,flato}, the geometric
quantization of quantum planes would show the difference in the roles
played be the parameter $q$ and the quantum Planck constant $\hbar$.

\end{document}